# CHAOS IN A SPATIAL EPIDEMIC MODEL[1]

By Rick Durrett and Daniel Remenik

*Cornell University*

We investigate an interacting particle system inspired by the gypsy moth, whose populations grow until they become sufficiently dense so that an epidemic reduces them to a low level. We consider this process on a random 3-regular graph and on the $d$-dimensional lattice and torus, with $d \geq 2$. On the finite graphs with global dispersal or with a dispersal radius that grows with the number of sites, we prove convergence to a dynamical system that is chaotic for some parameter values. We conjecture that on the infinite lattice with a fixed finite dispersal distance, distant parts of the lattice oscillate out of phase so there is a unique nontrivial stationary distribution.

**1. Introduction.** The inspiration for this paper arose almost 20 years ago. The first author had recently moved to Ithaca, New York and the Northeastern United States was in the midst of a gypsy moth infestation. For all of one summer, he and his wife destroyed egg masses, picked larvae off of trees and put bands of sticky tape to keep the larvae from climbing the trees. When the next summer came, the outlook for their trees seemed bleak, but suddenly all of the larvae were dead or deformed, a victim of the nuclear polyhedrosis virus, which spreads through the gypsy moth population once it becomes sufficiently dense.

To model this process we use dynamics that occur in discrete time with each site in some graph $G_N$ either occupied or vacant. The number of nodes in $G_N$ will be an increasing function of $N$ which tends to infinity. Two processes occur alternately: growth and epidemic.

*Growth.* Gypsy moths lay dormant in the winter as eggs, so no occupied site survives to the next time period but gives birth to a mean $\beta > 1$ number of individuals. Each individual born at $x$ is sent to a site randomly chosen

---

Received August 2008; revised November 2008.

[1]Supported in part by NSF Grant DMS-07-04996 from the probability program.

*AMS 2000 subject classifications.* 60K35, 60J10, 92D25, 37D45, 37N25.

*Key words and phrases.* Epidemic model, chaos, interacting particle system, dynamical system, random graph, gypsy moth.

from $\mathcal{N}_N(x) \subseteq G_N$, the *growth neighborhood of* $x$, which contains all of the nearest neighbors of $x$ in the graph but in general will be larger.

*Epidemic.* With a small probability $\alpha_N$ an infection lands at each site. If the site $x$ is occupied an infection starts which spreads from $x$ to all of its occupied neighbors in the graph and continues until all sites in the connected component of occupied sites containing $x$ are wiped out (observe that the larger the cluster of occupied sites, the more likely it is to be wiped out by the epidemic). It is assumed that the epidemic occurs rapidly so it is completed before the next growing season.

Our goal is to study this process on a random 3-regular graph and on a discrete torus of dimension $d \geq 2$. The second graph is more realistic from a biological point of view, but the first one is easier to deal with because explicit formulas are available. In both cases, infections will be transmitted along edges connecting neighbors. Observe that if we assume that $\alpha_N \to 0$ then only components with $O(1/\alpha_N)$ sites will be affected by epidemics. In site percolation on an regular tree of degree 3 and on $\mathbb{Z}^d$ there is phase transition from all components small to the existence of an infinite component at some density $p_c$. On the random 3-regular graph and the torus this phase transition produces one giant component of size $O(n)$. Thus we expect that the density of occupied sites will increase until $p > p_c$, at which point a large epidemic occurs and reduces the density to a low level and the cycle begins again. We will show that in some cases this leads to chaotic behavior of the densities.

### 1.1. *Mean-field growth on a random 3-regular graph.*

To work our way up to proving results about this system and the corresponding process on the torus we begin with the case in which $G_N$ is a random 3-regular graph with $N$ nodes, that is, a graph chosen at random from the set of graphs with $N$ vertices all of which have degree 3 ($N$ must be even). We will denote this random graph by $R_N$ and we will condition on the event that $R_N$ is connected. It is known [see Janson, Luczak and Rucinski (2000)] that the probability that $R_N$ is connected tends to 1. We choose this graph, not because it reflects reality, but because $R_N$ is locally a tree, so we have explicit formulas for the percolation probabilities. To have a simple process in which the number of occupied sites at the beginning of the growing season is a Markov process, we let $\mathcal{N}_N(x) = R_N$ for all $x$. As we will see, in the limit as $N \to \infty$ the result is a very interesting dynamical system.

To guess what this limiting system must be, observe that if we assume that the density of occupied sites before the growth step is $p$, so the expected number of occupied sites is $pN$, then the expected density after the birth step is

$$f_N(p) = 1 - \left(1 - \frac{\beta}{N}\right)^{pN} \approx f(p) = 1 - e^{-\beta p}.$$

Now the random 3-regular graph looks locally like a tree in which each vertex has degree 3 (we will refer to this tree as the 3-tree). Proceeding heuristically, in the limit $N \to \infty$, each occupied site survives the epidemic if and only if it is not in the giant component of the percolation process on the 3-tree defined by declaring open the sites that are occupied after the growth step. Thus if the density before the epidemic is $p$, the density $g_T(p)$ after the epidemic (the $T$ in the subscript is for tree) is exactly the probability that the origin is open in this percolation process but it does not percolate. The threshold for the existence of a giant component is $p_c = 1/2$, so if $p \leq 1/2$ then $g_T(p) = p$.

To compute the density for $p > 1/2$ we need to compute the percolation probability on the 3-tree. Throughout the rest of the paper, whenever we say percolation we mean the event that the origin is an infinite cluster of occupied sites. We start by noting that for site percolation on the binary tree (which is an infinite rooted tree where each vertex has two descendants, so all vertices have degree 3 except for the root which has degree 2) the percolation probability $\theta_{\text{bin}}(p)$ satisfies

$$\theta_{\text{bin}}(p) = p(1 - (1 - \theta_{\text{bin}}(p))^2)$$

since for this event to occur the origin must be occupied and percolation must occur from one of the two neighbors. Solving gives

$$\theta_{\text{bin}}(p) = \frac{2p - 1}{p} = 2 - \frac{1}{p}.$$

On the 3-tree the probability of percolation is then

$$\theta_T(p) = p(1 - (1 - \theta_{\text{bin}}(p))^3)$$

since the site must be occupied and percolation must occur from one of the three neighbors. Thus for $p \in (1/2, 1]$

$$g_T(p) = \mathbb{P}(0 \text{ is occupied}, |C_0| < \infty) = p - \theta_T(p) = p\left(\frac{1}{p} - 1\right)^3 = \frac{(1-p)^3}{p^2}.$$

Let $a_0$ be the solution of $1 - e^{-\beta a_0} = 1/2$ [i.e., $a_0 = (\log 2)/\beta$]. Combining the formulas for $f$ and $g_T$ we see that the limiting dynamical system should be the one defined by the function

$$h_T(p) = g_T(f(p)) = \begin{cases} 1 - e^{-\beta p}, & 0 \leq p \leq a_0, \\ \dfrac{e^{-3\beta p}}{(1 - e^{-\beta p})^2}, & a_0 < p \leq 1. \end{cases}$$

Observe that $h_T$ is continuous in $[0, 1]$.

We are interested in properties of the iterates of $h_T(p)$:

- If $\beta \leq 1$ then $f(p) < p$ for all $p > 0$ and thus $h_T^k(p)$ decreases to 0 as $k \to \infty$.

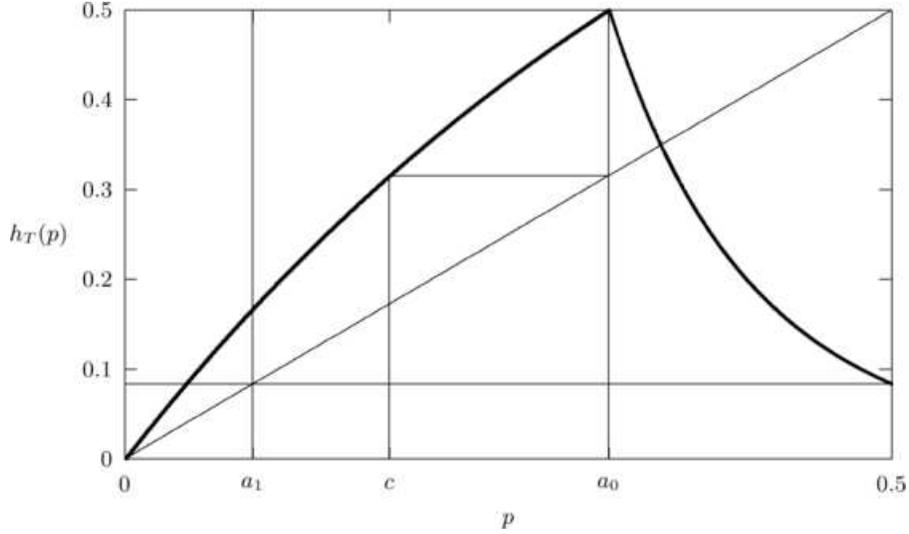

Fig. 1.  *Graph of $h_T$ with $\beta = 2\log 3$. The point $c = h_T^{-1}(a_0)$ will play a role in the proof of Theorem 1.*

- If $\beta > 1$ then starting from a small positive $p$, $f^k(p)$ increases to a unique fixed point $p^*$. If $p^* \le 1/2$ then we never get an epidemic and $h_T^k(p)$ increases to the same fixed point.
- $1/2$ is a fixed point when $e^{-\beta/2} = 1/2$, that is, $\beta = 2\log 2$. When $\beta > 2\log 2$, we let $a_1 = h_T(1/2) = e^{-3\beta/2}/(1 - e^{-\beta/2})^2$. Eventually the iterates of $h_T$ lie in the interval $[a_1, 1/2]$, and once they reach this interval, they stay there (see Figure 1).

Hence if $\beta \le \beta_c = 2\log 2$, $h_T(p) = f(p)$ for all $p$ and the epidemic part of the dynamics is not seen in the limiting system. If $\beta > \beta_c$ then $h_T(p) < 1/2 < f(p)$ for $p \ge a_0$.

Figure 2 shows the orbits of the system as a function of $\beta$. We plot $h_T^k(p)$ for $501 \le k \le 550$ to remove the initial transient. Note that the system proceeds directly from a stable fixed point to a "chaotic phase" rather than via period doubling bifurcations of the type occurring in the quadratic maps $rx(1-x)$. To say in what sense the behavior is chaotic, we will use two results of the theory of discrete time dynamical systems. The first result, which we include here for convenience, is commonly referred to as "period three implies chaos":

PROPOSITION 1.1 [Theorem 1 in Li and Yorke (1975)].  *Let $F: J \longrightarrow J$ be a continuous function on a real interval $J$ and assume that there is point $a \in J$ such that*

$$F^3(a) \le a < F(a) < F^2(a).$$

*Then:*

(a) *For every $k = 1, 2, \ldots$ there is a point in $J$ of period $k$, that is, a point $r \in J$ such that $F^k(r) = r$ but $F^j(r) \neq r$ for $0 < j < k$.*

(b) *There is an uncountable set $S \subseteq J$ containing no periodic points such that:*

    (b.i) *For every $p, q \in S$, $p \neq q$,*

$$\limsup_{N \to \infty} |F^N(p) - F^N(q)| > 0$$

    *and*

$$\liminf_{N \to \infty} |F^N(p) - F^N(q)| = 0.$$

    (b.ii) *For every $p \in S$ and any periodic point $q \in J$,*

$$\limsup_{N \to \infty} |F^N(p) - F^N(q)| > 0.$$

We will say that $F$ is *chaotic* if $F$ satisfies the conditions (a) and (b) above. (b.ii) rules out convergence to periodic orbits, while (b.i) shows that all the points in $S$ have different limiting behaviors.

### Theorem 1.

(a) *The dynamical system defined by the function $h_T : [a_1, 1/2] \longrightarrow [a_1, 1/2]$ is chaotic for every $\beta > 2 \log 2$.*

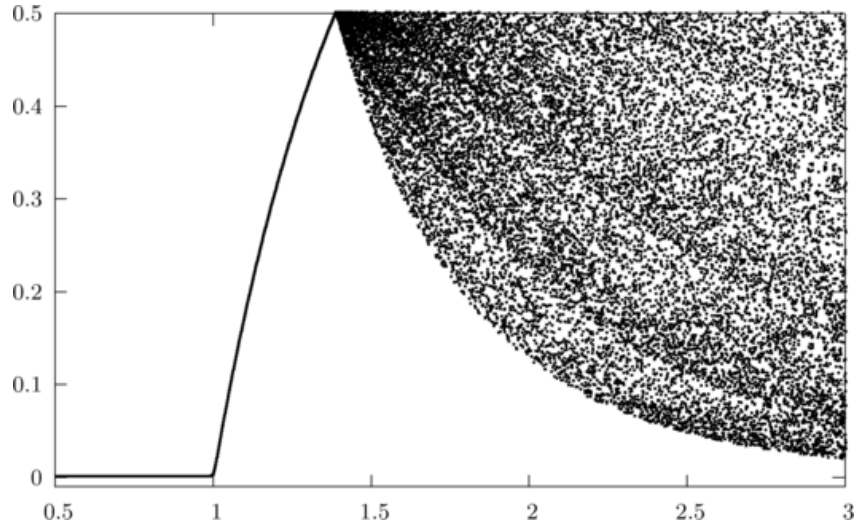

Fig. 2. *Orbits of the system $(h_T^k(p))_{k \geq 0}$ started at $p = 0.1$. The x-axis has the values of $\beta$ used in the simulations, while the y-axis has $h_T^k(p)$ for $k = 501, \ldots, 550$.*

(b) *If $\beta \in (2\log 2, 2.48]$ then the system has an invariant measure, $\mu = \mu \circ h_T^{-1}$, which is absolutely continuous with respect to the Lebesgue measure.*

Simulations suggest that (b) actually holds for all $\beta > 2\log(2)$.

Now we come back to the process running on $R_N$. We will denote our process by $\eta_k^N$, with $\eta_k^N(i) = 1$ if $i$ is occupied at time $k$ and $\eta_k^N(i) = 0$ if not. The density of occupied sites at time $k$ will be denoted by $\rho_k^N$:

$$(1.1) \qquad \rho_k^N = \frac{1}{N}|\eta_k^N| = \frac{1}{N}\sum_{i=1}^N \eta_k^N(i).$$

The initial distribution $\eta_0^N$ of the process will always be assumed to be a product measure with some density $p \in [0,1]$ (so, in particular, $\rho_0^N$ converges in probability to $p$). In the preceding discussion we argued heuristically that $\rho_k^N$ converges to the deterministic system defined by $h_T$. The next result shows that this is indeed the case:

THEOREM 2. *Assume that $G_N = R_N$ and that the infection probability of the epidemic satisfies*

$$\alpha_N \log_2 N \xrightarrow[N\to\infty]{} \infty.$$

*Then the process $(\rho_k^N)_{k\geq 0}$ converges in distribution to the (deterministic) orbit, starting at $p$, of the dynamical system associated to $h_T$.*

The above convergence means that $(\rho_k^N)_{k\geq 0}$ converges in distribution to a deterministic process whose paths are given by the orbits $(h_T^k(p))_{k\geq 0}$.

1.2. *Local growth on the $d$-dimensional torus.* Turning now to a more realistic setting, we consider the process running on the $d$-dimensional torus $(\mathbb{Z} \bmod N)^d$, for $d \geq 2$, which we will denote by $\mathbb{T}_N$. The case $d = 2$ is the one relevant to gypsy moths, but it is no harder to prove our results in general.

To add some more realism and make our process more interesting, we will take now the growth neighborhoods $\mathcal{N}_N(x)$ to be smaller than $\mathbb{T}_N$. We let

$$\mathcal{N}_N(x) = \{y \in \mathbb{T}_N : 0 < \|y - x\|_\infty \leq r_N\}$$

(here the difference $y - x$ is computed modulo $N$) and take the range $r_N$ to be such that $r_N \to \infty$. (We remark that on $\mathbb{T}_N$ we are considering the $L^1$ distance; in particular, two points $x, y \in \mathbb{T}_N$ are neighbors if $\|x - y\|_1 = 1$).

We start as before by guessing what the limiting system should be. To do this we will assume for a moment that $r_N = \infty$ for all $N$, so we are back in the case of mean-field growth of the previous subsection. The growth step behaves exactly as before: if $p$ is the density of occupied sites before the growth step, then the density after is

$$f_{N^2}(p) = 1 - \left(1 - \frac{\beta}{N^2}\right)^{pN^2} \approx f(p) = 1 - e^{-\beta p}.$$

The behavior of the epidemic step in the limit $N \to \infty$ is analogous to the one in the random 3-regular graph: if $p$ is the density of occupied sites before the epidemic, then the density $g_L(p)$ after (here the subscript $L$ is for lattice) is the probability that the origin is open but does not percolate in a site percolation process in $\mathbb{Z}^d$.

Unlike the case of percolation on the 3-tree, we do not have an explicit formula available for the percolation probability in $\mathbb{Z}^d$, but we still know some qualitative properties. Letting $C_0$ be the percolation cluster containing the origin and

$$\theta_L(p) = \mathbb{P}(|C_0| = \infty)$$

we have that there is a $p_c \in (0,1)$ ($p_c \approx 0.593$ in $d = 2$) such that $\theta_L(p) = 0$ for $p \le p_c$, $\theta_L(p)$ is strictly increasing for $p > p_c$, and $\theta_L(p)$ is infinitely differentiable at every $p \ne p_c$ [see Theorem 8.92 of Grimmett (1999)]. We also have that

$$g_L(p) = \mathbb{P}(0 < |C_0| < \infty) = \mathbb{P}(|C_0| < \infty) - \mathbb{P}(|C_0| = 0) = p - \theta_L(p),$$

so $g_L(p)$ is infinitely differentiable at $p \ne p_c$ and $g_L(p) = p$ for $p \le p_c$.

As before we let $h_L(p) = g_L(f(p))$ and $\beta_c$ be the value of $\beta$ solving $p_c = 1 - e^{-\beta p_c}$, that is,

$$\beta_c = \frac{1}{p_c} \log\left(\frac{1}{1 - p_c}\right)$$

($\beta_c \approx 1.516$ in $d = 2$). Observe that $g_L(p) \in (0,1)$ for $p \in (0,1)$ so, in particular, $h_L(p) > 0$ for $p > 0$. Our next result holds under an hypothesis on the percolation function which might seem strange at a first look, but which holds in $d = 2$ and is expected to also hold in $3 \le d < 6$.

Theorem 3. *Suppose that*

$$(1.2) \qquad \lim_{p \downarrow p_c} \theta_L'(p) = \infty.$$

*Then there is an $\varepsilon > 0$ such that for every $\beta \in (\beta_c, \beta_c + \varepsilon)$ the dynamical system $(h_L^k(p))_{k \ge 0}$ has an invariant measure which is absolutely continuous with respect to the Lebesgue measure.*

We believe (and simulations suggest) that the result holds for all $\beta > \beta_c$. As Yuval Peres pointed out to us, it is easy to show that (1.2) holds in $d = 2$ using Russo's formula and the fact that the expected number of pivotal sites goes to infinity as $p \downarrow p_c$ in two dimensions. This argument would obviously work in other dimensions too if we knew that the expected number of pivotal sites blows up at $p_c$. This should be the case in $3 \leq d < 6$ because it is expected that $\theta_L(p) \approx C(p - p_c)^\gamma$ as $p \downarrow p_c$ with $\gamma < 1$ in $d < 6$, $\gamma = 1$ in $d > 6$, and with logarithmic corrections in $d = 6$ [see, e.g., Chapter 9 of Grimmett (1999)].

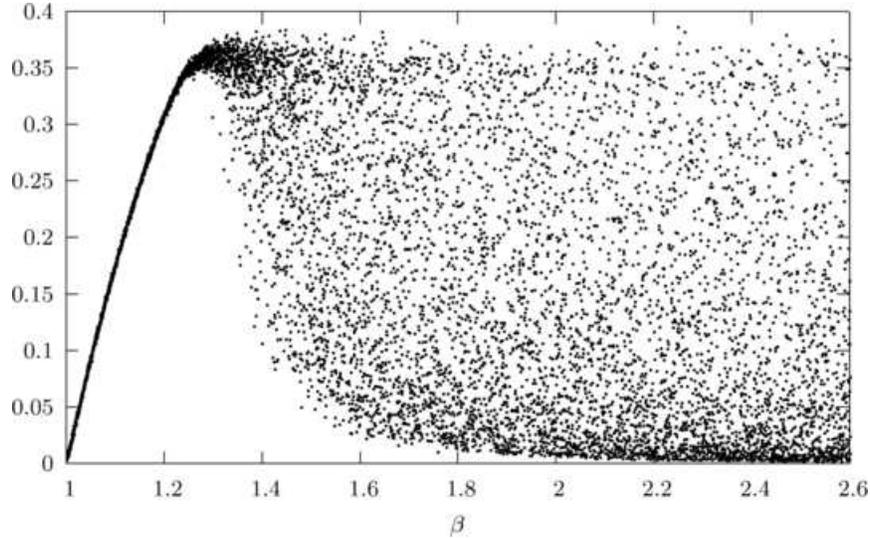

Fig. 3. *Density process* $(\rho_k^N(p))_{k \geq 0}$ *running on the two-dimensional torus with mean-field growth, depicted for* $k = 501, \ldots, 550$ *for different values of the parameter* $\beta$ *(similar to Figure 2). Here* $N = 500$.

Our next goal is to show that the process $\rho_k^N$ on the torus $\mathbb{T}_N$ converges to the deterministic orbit of the dynamical system defined by $h_L$. The processes $\eta_k^N$ and $\rho_k^N$ are defined in this case exactly as for the random 3-regular graph; see (1.1) and the preceding lines. If we consider the case of mean-field growth [i.e., $\mathcal{N}_N(x) = \mathbb{T}_N$] then the result follows from the same arguments as those we will use to prove Theorem 2 (the proof is actually simpler because we do not have to prove that the torus looks locally like $\mathbb{Z}^d$). Figure 3 shows part of the trajectories of $\rho_k^N$ in the case of mean-field growth. But, as we mentioned, we want to deal with the more general case $\mathcal{N}_N(x) = \{y \in \mathbb{T}_N : 0 < \|y - x\|_\infty \leq r_N\}$ with $r_N \to \infty$. The result does not seem to be true if we do not take $r_N \to \infty$. As Figure 4 shows, the graph of $\{(\rho_k^N, \rho_{k+1}^N), k \geq 0\}$ does not correspond to any function. This difficulty disappears as $N \to \infty$ if we take $r_N \to \infty$ at an appropriate rate.

We will assume the following on $\alpha_N$ and $r_N$:

$$\frac{r_N}{N} \longrightarrow 0$$

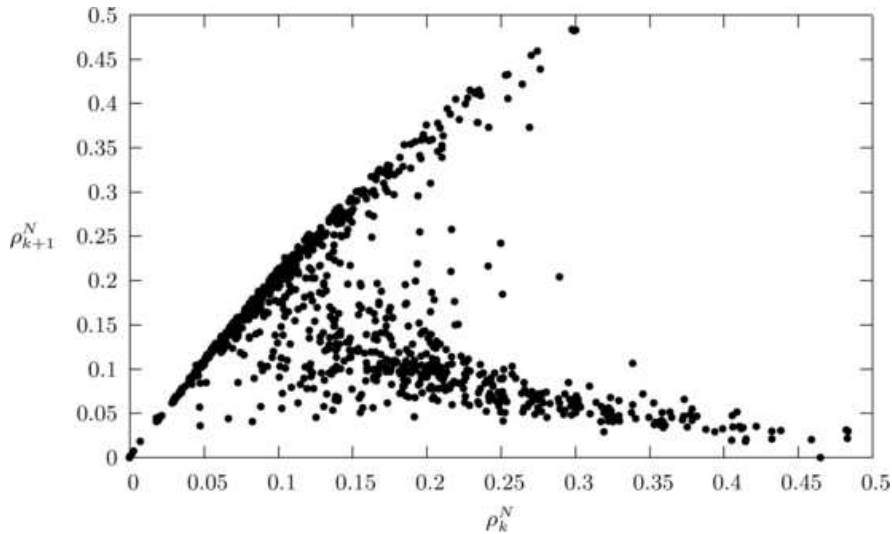

Fig. 4. *Graph of $\rho_k^N$ against $\rho_{k+1}^N$ on the two-dimensional torus with $N = 750$ and $r_N = 50$. The graph clearly does not correspond to a function.*

and

$$\alpha_N r_N \longrightarrow \infty.$$

For instance, we could take $r_N = N^\gamma$ and $\alpha_N = N^{-\delta}$ for some $0 < \delta < \gamma < 1$.

THEOREM 4. *Assume that $G_N = \mathbb{T}_N$, with $d \geq 2$, and that the number of individuals to which each occupied site gives birth to during the growing season is a Poisson random variable with mean $\beta$. Then the process $(\rho_k^N)_{k \geq 0}$ converges in distribution to the (deterministic) orbit, starting at $p$, of the dynamical system associated to $h_L$.*

1.3. *Local growth on $\mathbb{Z}^d$.* We now consider the case in which $r_N$ is constant. Figure 5 shows that when $r_N = 5$ the fluctuations in the density of occupied sites decrease as the system size increases. Figure 6 shows a picture of the process running on the torus of size $450 \times 450$ with $r_N = 5$. As this picture suggests the density stays constant because different parts of the lattice oscillate out of phase.

THEOREM 5. *Consider the process running in $\mathbb{Z}^d$ with $d \geq 2$. If $r_N = L$ and $L$ is sufficiently large then there is a nontrivial stationary distribution.*

SKETCH OF THE PROOF. The key to the proof is that the density of occupied sites after growth is at most $f(1) = 1 - e^{-\beta}$ so after the epidemic there will be a positive density of occupied sites. Let $\delta = (1 - e^{-\beta})e^{-4\beta}$ be the probability that a site is occupied and has four vacant neighbors. Divide space into squares of side $L/2$ and declare that the square is occupied if at least a fraction $\delta/2$ of the sites are. If $L$ is large enough and $T$ is chosen suitably then the set of occupied squares at time $nT$ dominates oriented percolation with $p$ close to 1 and the result follows from standard "block construction" arguments [for an account of this method see, e.g., Durrett (1995)]. By order of the Associate Editor further details are left to the reader. □

The remainder of the paper is devoted to proofs. The proof of Theorem 1 is given in Section 2. If you get bored with all of the algebra and calculus involved you can skip to Section 3 where the proof of Theorem 2 is given. The proof of Theorem 3 given in Section 4 and the more complicated proof of Theorem 4 in Section 5 rely on ideas from Sections 2 and 3, but are independent of each other.

The authors would like to thank referee Nicolas Lanchier for his careful reading of the paper which resulted in a number of corrections and clarifications.

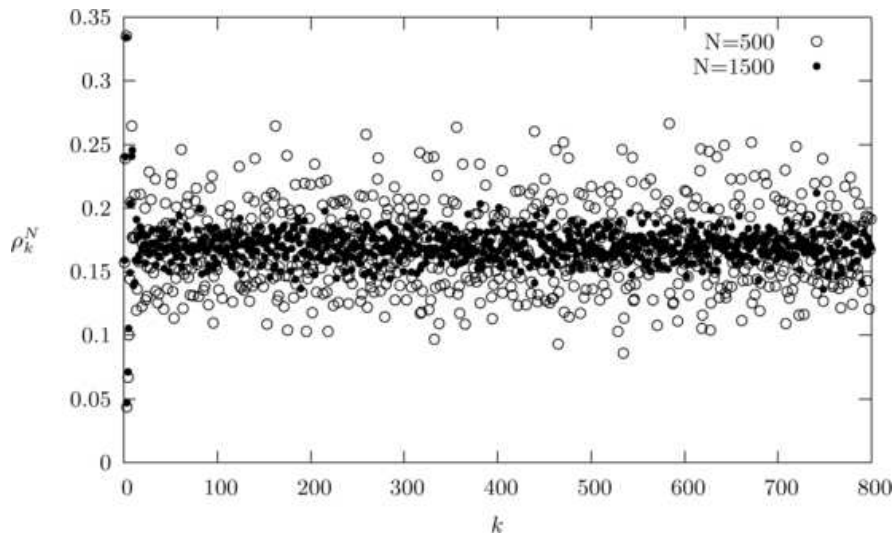

FIG. 5. Sequence of densities $\rho_k^N$ of the process running on the two-dimensional torus with local interactions in the epidemic step for $N = 500$ and $N = 1500$, both with $r_N = 5$. As this graph suggests, the fluctuations of the density process get small as $N$ grows if the range $r_N$ is held fixed.

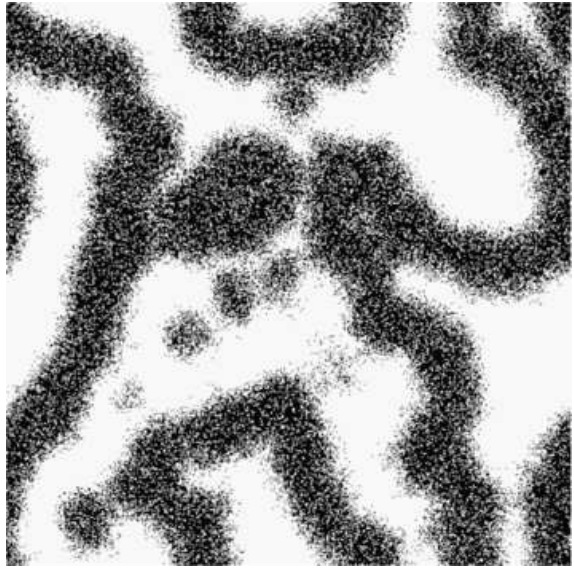

Fig. 6. *State of the process at time 200 on a torus of size $450 \times 450$ (black dots are occupied). In this simulation, $\beta = 2.25$, $r_N = 5$ and the infection probability at each site is $5 \cdot 10^{-6}$. This picture corresponds to an intermediate state of the process, after an epidemic event wiped out a big cluster but the process has had time to grow back.*

**2. Proof of Theorem 1.** By Proposition 1.1, to obtain (a) it is enough to prove that there is a point $c \in [a_1, 1/2]$ such that

$$h_T^3(c) \le c < h_T(c) < h_T^2(c).$$

In our case we can take

$$c = f^{-1}(a_0) = \frac{1}{\beta} \log\left(\frac{\beta}{\beta - \log 2}\right)$$

(see Figure 1). Observe that since $a_0 < 1/2$, $c = \beta^{-1} \log((1-a_0)^{-1}) < \beta^{-1} \times \log 2 = a_0$. Hence

$$h_T(c) = f(c) = a_0,$$
$$h_T^2(c) = f(a_0) = \tfrac{1}{2}$$

and

$$h_T^3(c) = h_T(1/2) = a_1.$$

It is clear then that $c < h_T(c) < h_T^2(c)$. To see that $h_T^3(c) \le c$ we need to show that $a_1 \le f^{-1}(a_0)$, that is, that

$$\frac{e^{-3\beta/2}}{(1 - e^{-\beta/2})^2} \le \frac{1}{\beta} \log\left(\frac{\beta}{\beta - \log 2}\right)$$

or, equivalently, that

$$(2.1) \qquad \phi_1(\beta) = \exp\left(\frac{\beta e^{-3\beta/2}}{(1 - e^{-\beta/2})^2}\right) \leq \phi_2(\beta) = \frac{\beta}{\beta - \log 2}$$

for all $\beta > 2 \log 2$. If you look at the picture of these two functions it seems clear that the inequality holds, but the proof is not as simple as the picture suggests. We will divide it into two parts.

First, assume that $\beta \in (2 \log 2, 1.75]$. We will show that

$$(2.2) \qquad \phi_1(\beta) \leq 4 - \frac{\beta}{\log 2} \leq \phi_2(\beta).$$

To get the first inequality let

$$\sigma(\beta) = \frac{\beta e^{-3\beta/2}}{(1 - e^{-\beta/2})^2}.$$

A simple calculation gives

$$\sigma''(\beta) = \frac{9e^\beta - 4e^{\beta/2} + 1}{4e^{5\beta/2} - 16e^{2\beta} + 24e^{3\beta/2} - 16e^\beta + 4e^{\beta/2}}$$

and we claim that this quotient is positive. Indeed, it is easy to see that the numerator is positive, while putting $a = e^{\beta/2}$ the denominator becomes $4a^5 - 16a^4 + 24a^3 - 16a^2 + 4a$, so dividing by $4a$ we need to show that

$$w(a) = a^4 - 4a^3 + 6a^2 - 4a + 1 > 0$$

for all $a > 2$. Observe that $w'(a) = 4a^3 - 12a^2 + 12a - 4$, so $w'(2) = 4$, while $w''(a) = 12(a-1)^2 > 0$, so $w'(a) > 0$ for all $a > 2$. Since $w(2) = 1$ we deduce that $w(a) > 0$ for all $a > 2$ as required. Hence $\sigma$ is convex, and thus so is $\phi_1 = \exp(\sigma(\cdot))$. Since

$$\phi_1(2\log(2)) = 2 = 4 - \frac{2\log 2}{\log 2} \quad \text{and} \quad \phi_1(1.75) \approx 1.4518 < 4 - \frac{1.75}{\log 2} \approx 1.4753,$$

the convexity of $\phi_1$ gives the desired inequality.

To get the second inequality in (2.2), observe that

$$\phi_2(2\log 2) = 2 = 4 - \frac{2\log 2}{\log 2} \quad \text{and} \quad \phi_2'(2\log 2) = -\frac{1}{\log 2}.$$

Therefore, since this last quantity is exactly the slope of the line appearing in the middle term of (2.2) and since $\phi_2$ is strictly convex, we deduce that $\phi_2'(\beta)$ is larger than this slope for every $\beta > 2 \log 2$ and thus the inequality holds.

Now we assume that $\beta > 1.75$. Using the Taylor expansion of the functions $1/(1-x)$ and $e^x$ about $x = 0$ we get that (2.1) is equivalent to

$$\sum_{n \geq 0} \left(\frac{\log 2}{\beta}\right)^n \geq \sum_{n \geq 0} \frac{1}{n!}\left(\frac{\beta e^{-3\beta/2}}{(1 - e^{-\beta/2})^2}\right)^n,$$

so it is enough to show that

$$\left(\frac{\log 2}{\beta}\right)^n \geq \frac{1}{n!}\left(\frac{\beta e^{-3\beta/2}}{(1-e^{-\beta/2})^2}\right)^n$$

for all $n \geq 0$ and $\beta > 1.75$. The inequality holds trivially for $n = 0$, so by induction it is enough to prove that

$$\frac{\log 2}{\beta} \geq \frac{1}{n} \frac{\beta e^{-3\beta/2}}{(1-e^{-\beta/2})^2}$$

for all $n \geq 1$ or, equivalently, for $n = 1$. That is, we need to show that

$$(2.3) \qquad \beta^2 \frac{e^{-3\beta/2}}{(1-e^{-\beta/2})^2} \leq \log 2$$

for all $\beta > 1.75$. To see that this holds we observe that the derivative of the left side with respect to $\beta$ is

$$-\frac{\beta e^{-\beta/2}(3\beta e^{\beta/2} - 4e^{\beta/2} - \beta + 4)}{2(e^{\beta/2} - 1)^3}.$$

We claim that this quotient is negative for $\beta > 1.75$. Indeed, the denominator is clearly positive, so we only need to show that

$$w(\beta) = 3\beta e^{\beta/2} - 4e^{\beta/2} - \beta + 2 > 0$$

for $\beta > 1.75$. This is easy, because $w'(\beta) = 3e^{\beta/2}(1+\beta/2) - 2e^{\beta/2} - 1 > e^{\beta/2} - 1 > 0$ and $w(1.75) \approx 5.28$. Thus the left side of the (2.3) is decreasing in $\beta$, and then the inequality holds because its value at $\beta = 1.75$ is approximately $0.6523 < \log 2$. This finishes the proof of (a).

To get (b) it is enough to show by Lasota and Yorke (1973) that

$$(2.4) \qquad \inf_{p\in[a_1,1/2]\setminus\{a_0\}} |(h_T^3)'(p)| > 1$$

for $\beta \in (2\log 2, 2.48]$. The idea of the proof is the following. We find an explicit formula for $(h_T^3)'$ and use it to compute numerically its infimum on $[a_1, 1/2] \setminus \{a_0\}$ for every $\beta$ in a certain grid of $(2\log 2, 2.48]$. Due to monotonicity properties of the derivative of $h_T$ the numerical computation of the infimum is exact (up to floating-point numerical errors which are small enough for our purposes) for any fixed $\beta$. We then show that $(h_T^3)'$, as a function of $\beta$, has a Lipschitz constant that ensures that the infimum is larger than 1 for every $\beta$ between subsequent points in the grid. We will do this step by step.

We begin by computing $(h_T^3)'$. For $p \in [a_1, a_0)$, $h_T'(p) = f'(p) = \beta e^{-\beta p}$, while for $p \in (a_0, 1/2]$,

$$h_T'(p) = \frac{-3\beta e^{-3\beta p}}{(1-e^{-\beta p})^2} - 2\frac{e^{-3\beta p}}{(1-e^{-\beta p})^3}\beta e^{-\beta p} = \frac{e^{-3\beta p}}{(1-e^{-\beta p})^3}[-3\beta + \beta e^{-\beta p}].$$

This gives an explicit formula for $h'_T$. On the other hand,

$$(2.5) \qquad (h^3_T)'(p) = h'_T(h^2_T(p))h'_T(h_T(p))h'_T(p).$$

Putting these two formulas together we get an explicit expression for $(h^3_T)'$.

Now observe that $h'_T$ is decreasing in $[a_1, a_0)$ and increasing in $(a_0, 1/2]$. Indeed, $h''_T(p) = f''(p) = -\beta^2 e^{-\beta p} < 0$ on the first interval, while on the second one $h''_T(p) = g''_T(f(p))f'(p)^2 + g'_T(f(p))f''(p)$, so since $f' > 0$, $f'' < 0$,

$$g'_T(p) = \left(\frac{1}{p} - 1\right)^3 + 3p\left(\frac{1}{p} - 1\right)^2\left(\frac{-1}{p^2}\right) = -\left(1 + \frac{2}{p}\right)\left(\frac{1}{p} - 1\right)^2 < 0$$

and

$$g''_T(p) = \frac{2}{p^2}\left(\frac{1}{p} - 1\right)^2 - \left(1 + \frac{2}{p}\right)\left(\frac{1}{p} - 1\right)\left(\frac{-1}{p^2}\right) > 0,$$

we get that $h''_T(p) > 0$ for $p \in (a_0, 1/2]$. This means by (2.5) that $(h^3_T)'$ is monotone in each interval of constancy of its sign. These intervals are given by the partition of $[a_1, 1/2]$ defined by the preimage of $a_0$ under $h^3_T$. We deduce that

$$\inf_{p \in [a_1, 1/2]\backslash\{a_0\}} |(h^3_T)'(p)| = \inf_{p \in h^{-3}_T(a_0)\cup\{a_1, 1/2\}} \min\{|(h^3_T)'(p^-)|, |(h^3_T)'(p^+)|\},$$

where the superscripts $-$ and $+$ indicate left and right derivatives, respectively. Using this observation we can compute numerically the infimum in (2.4) for any given $\beta$. We did this for every $\beta$ in a grid of width $2 \cdot 10^{-6}$ of $(2\log 2, 2.48]$, and we obtained that the infimum is larger than 1.002 at each of these values of $\beta$. Figure 7 shows a graph of the values obtained.

The last step is to make sure that the infimum in (2.4) stays above 1 for every $\beta \in (2\log 2, 2.48]$. We will write $h_T(p, \beta)$ to indicate the dependence of $h_T(p)$ on the value of the parameter $\beta$. Our goal is to find a bound for $|\frac{\partial^2}{\partial\beta\,\partial p}h^3_T(p, \beta)|$. Observe that by the product rule and (2.5), if $|\frac{\partial}{\partial p}h_T(p, \beta)| \leq M_1$ and $|\frac{\partial^2}{\partial\beta\,\partial p}h_T(p, \beta)| \leq M_2$ for all $\beta \in (2\log 2, 2.48]$ and $p \in [a_1, 1/2] \setminus \{a_0\}$ then

$$(2.6) \qquad \left|\frac{\partial^2}{\partial\beta\,\partial p}h^3_T(p, \beta)\right| \leq 3M_1^2M_2$$

for all such $\beta$ and $p$. We already computed $|\frac{\partial}{\partial p}h_T(p, \beta)|$. For $p \in [a_1, a_0)$, it equals $\beta e^{-\beta p}$ which is smaller than 2.48 for each $\beta \leq 2.48$. For $p \in (a_0, 1/2]$ we know that $h'_T$ is negative and increasing, so

$$\left|\frac{\partial}{\partial p}h_T(p, \beta)\right| \leq \left|\frac{\partial}{\partial p}h_T\left(\frac{1}{2}, \beta\right)\right| = \left|\frac{e^{-3\beta/2}}{(1 - e^{-\beta/2})^3}[-3\beta + \beta e^{-\beta/2}]\right|$$

$$\leq \frac{e^{-3\cdot 2.48/2}}{2^{-3}} \cdot 4 \cdot 2.48 \approx 1.923.$$

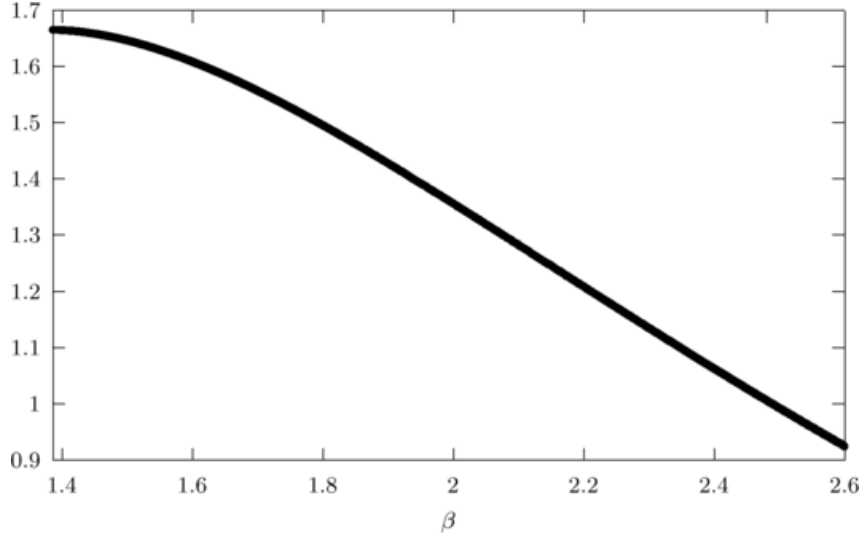

Fig. 7. *Infimum of* $|(h_T^3)'(p)|$ *on the relevant interval for* $\beta \in (2\log 2, 2.6)$. *The computation was done for each* $\beta$ *on a grid of width* $2 \cdot 10^{-6}$ *on this interval, as explained within the proof of Theorem 1. The infimum lies above 1.002 for* $\beta \in (2\log 2, 2.48]$.

Thus if we take $M_1 = 2.48$ the desired inequality holds. Now for $p \in [a_1, a_0)$,

$$\left| \frac{\partial^2}{\partial \beta \, \partial p} h_T(p, \beta) \right| = \left| \frac{\partial}{\partial \beta} (\beta e^{-\beta p}) \right| = |(1 - \beta^2) e^{-\beta p}| \leq 1.$$

For $p \in (a_0, 1/2]$,

$$\left| \frac{\partial^2}{\partial \beta \, \partial p} h_T(p, \beta) \right|$$

$$= \left| \frac{\partial}{\partial \beta} \left( \frac{e^{-3\beta p}}{(1 - e^{-\beta p})^3} [-3\beta + \beta e^{-\beta p}] \right) \right|$$

$$= \frac{e^{-\beta p}}{(1 - e^{-\beta p})^4} |(9\beta p - 3) e^{-2\beta p} + (4 - 4\beta p) e^{-3\beta p} + (\beta p - 1) e^{-4\beta p}|$$

$$\leq \frac{e^{-\beta a_0}}{(1 - e^{-\beta/2})^4} (14\beta p + 8) \leq \frac{2^{-1}}{(1 - e^{-2.48/2})^4} (14 \cdot 2.48/2 + 8) \approx 49.73,$$

if $\beta \in (2\log 2, 2.48]$. Thus if we take $M_2 = 49.73$ we get by (2.6) that $|\frac{\partial}{\partial \beta} h_T^3(p, \beta)| \leq 917.6$.

The bound we just obtained implies that for any fixed $p \in [a_1, 1/2] \setminus \{a_0\}$ the function $\beta \mapsto \frac{\partial}{\partial p} h_T(p, \beta)$ is Lipschitz and its Lipschitz constant is at most 917.6. Now fix $\beta \in (2\log 2, 2.48]$ and let $\beta'$ be the point in the grid of $(2\log 2, 2.48]$ on which we computed the infimum in (2.4) which is immediately before $\beta$. Then for any $p \in [a_1, 1/2] \setminus \{a_0\}$,

$$\left|\frac{\partial}{\partial p} h_T^3(p, \beta)\right| \geq \left|\frac{\partial}{\partial p} h_T^3(p, \beta')\right| - \left|\frac{\partial}{\partial p} h_T^3(p, \beta) - \frac{\partial}{\partial p} h_T^3(p, \beta')\right|$$

$$\geq 1.002 - 917.6|\beta - \beta'| \geq 1.002 - 917.6 \cdot 2 \cdot 10^{-6} \approx 1.0001.$$

This completes the proof of (2.4).

**3. Proof of Theorem 2.** To prove this result it will be enough to study the one-step transition probabilities for $\rho_k^N$. Recall that in the growth step, since here $\mathcal{N}_N(x) = G_N$, every site becomes occupied with probability $1 - (1 - \beta/N)^{pN} \approx 1 - e^{-\beta p}$, where $p$ is the starting density of occupied sites. For simplicity we will assume that the occupation probability of each site after the growth step is exactly $1 - e^{-\beta p}$, and then in the proof of the theorem we will say how to remove this assumption.

Abusing notation, we will also let $\eta_k^N$ stand for the set of occupied sites in the process. $\eta_{k+1/2}^N$ will denote the intermediate state of the process between $\eta_k^N$ and $\eta_{k+1}^N$ after the growth part of the dynamics has been run but before running the epidemic. We will denote by $\{0, \ldots, N-1\}$ the set of nodes of $R_N$. $B(i, r)$ will denote the set of sites in $R_N$ at distance at most $r$ from $i$ (here the distance between two points $i$ and $j$ is defined as the number of edges in the shortest path going from $i$ to $j$).

Let $\tilde{\eta}_1^N$ be the set of occupied sites after the epidemic is run on $\eta_{1/2}^N$ ignoring infections coming from a distance greater that $(\log_2 N)/5$. Define $\tilde{\rho}_1^N = |\tilde{\eta}_1^N|/N$. Recall that we are assuming that

$$\alpha_N \log_2 N \longrightarrow \infty.$$

LEMMA 3.1.

$$\mathbb{E}(|\tilde{\rho}_1^N - \rho_1^N|) \xrightarrow[N \to \infty]{} 0,$$

uniformly in the initial density $p$.

PROOF. By translation invariance, and observing that $\tilde{\eta}_1^N(i) \geq \eta_1^N(i)$ for all $i \in R_N$,

$$\mathbb{E}(|\tilde{\rho}_1^N - \rho_1^N|) \leq \frac{1}{N} \sum_{i \in R_N} \mathbb{E}(|\tilde{\eta}_1^N(i) - \eta_1^N(i)|) = \mathbb{P}(0 \in \tilde{\eta}_1^N) - \mathbb{P}(0 \in \eta_1^N)$$

$$= \mathbb{P}(0 \in \tilde{\eta}_1^N \setminus \eta_1^N) \leq (1 - \alpha_N)^{1/5 \log_2 N} \approx e^{-1/5 \alpha_N \log_2 N} \longrightarrow 0.$$

The second inequality above follows from the fact that if 0 is in $\tilde{\eta}_1^N$ but not in $\eta_1^N$, then there must be an open path in $\eta_{1/2}^N$ going from 0 to $\partial B(0, (\log_2 N)/5)$, and all sites in this path must have not been infected. □

Now let

$$H_N = \{i \in R_N : B(i, (\log_2 N)/5) \text{ is a finite 3-tree}\}.$$

By a finite 3-tree we mean a finite tree where all nodes have degree 3 except for the leaves which have degree 1. The next lemma says that $R_N$ looks locally like a 3-tree:

LEMMA 3.2.

$$\mathbb{E}\left(\frac{1}{N}|R_N \setminus H_N|\right) = \mathbb{P}(0 \notin H_N) \xrightarrow[N \to \infty]{} 0.$$

PROOF. A random 3-regular graph is a special case of a graph with a fixed degree distribution and can be studied using techniques in Section 3.2 of Durrett (2007). To explore the subgraph $B(0, (\log_2 N)/5)$ of $R_N$, let $R_0 = \varnothing$, $A_0 = \{0\}$ and $U_0 = \{1, \ldots, N-1\}$. These are called the removed, active and unexplored sites, respectively. If $A_n \neq \varnothing$ then to go from time $n$ to $n+1$ we pick a site $i_n$ from $A_n$ according to some given rule and let

$$R_{n+1} = R_n \cup \{i_n\},$$

$$A_{n+1} = (A_n \setminus \{i_n\}) \cup \{j \in U_n : j \sim i\},$$

$$U_{n+1} = U_n \setminus \{j \in U_n : j \sim i\},$$

where $j \sim i$ here denotes that $j$ and $i$ are neighbors. For $n \leq 3N^{1/5}/2$, $|A_n| \leq 3N^{1/5}/2 + 2$, so the probability of a collision (i.e., that when we examine the neighbors of $i_n$ we see a site already in $A_n$) at some time is at most

$$2 \cdot \frac{3}{2}N^{1/5}\frac{3N^{1/5}/2 + 2}{N} \longrightarrow 0.$$

Now suppose that when choosing the sites $i_n$ we choose those at distance 1 from 0 first, then those at distance 2, etcetera. Then by time $3N^{1/5}/2$ we will have investigated all points within distance $(\log_2 N)/5$ of 0, and if we see no collision, then we will know that the subgraph $B(0, (\log_2 N)/5)$ is a tree. □

LEMMA 3.3. *Let $C_0$ be the cluster containing the origin in a site percolation process on the 3-tree, and let $\mathbb{P}_p$ denote the law of this process when each site is retained independently with probability $p \in [0, 1]$. Then for any $k_N \uparrow \infty$,*

$$\sup_{p \in [0,1]} |\mathbb{P}_p(\mathrm{diam}(C_0) < \infty) - \mathbb{P}_p(\mathrm{diam}(C_0) \leq k_N)| \xrightarrow[N \to \infty]{} 0.$$

PROOF. The result follows from the fact that any increasing sequence of continuous functions on $[0, 1]$ which converges pointwise to a continuous function on $[0, 1]$ actually converges uniformly to that function [see, e.g., Theorem 7.13 in Rudin (1976)]. We only need to observe that $\mathbb{P}_p(\mathrm{diam}(C_0) < \infty)$ and $\mathbb{P}_p(\mathrm{diam}(C_0) \leq k_N)$ are continuous on $[0, 1]$ as functions of $p$, and the latter is increasing in $N$ and converges pointwise to the former as $N \to \infty$. $\square$

LEMMA 3.4.

$$\mathbb{E}\left(\frac{1}{N}|\widetilde{\eta}_1^N \cap H_N|\right) \xrightarrow[N \to \infty]{} h_T(p),$$

uniformly in the initial density $p$.

PROOF. Observe that since $0 \in \widetilde{\eta}_1^N$ implies that $0 \in \widetilde{\eta}_{1/2}^N = \eta_{1/2}^N$,

$$
\begin{aligned}
(3.1) \quad & \mathbb{E}\left(\frac{1}{N}|\widetilde{\eta}_1^N \cap H_N|\right) \\
& \qquad = \mathbb{P}(0 \in \widetilde{\eta}_1^N \mid 0 \in H_N \cap \eta_{1/2}^N)\mathbb{P}(0 \in H_N)\mathbb{P}(0 \in \eta_{1/2}^N).
\end{aligned}
$$

By Lemma 3.2, $\mathbb{P}(0 \in H_N) \to 1$ uniformly in $p$, while by our assumption, $\mathbb{P}(0 \in \eta_{1/2}^N) = 1 - e^{-\beta p}$.

For the other term on the right-hand side of (3.1), we only need to look at the configuration of $\eta_{1/2}^N$ inside $B(0, (\log_2 N)/5)$, on which, conditional on the event $\{0 \in H_N\}$, the graph looks like a finite 3-tree. Thus, we can construct the random variables $(\widetilde{\eta}_1^N(0))_{N>0}$ conditioned on $\{0 \in H_N \cap \eta_{1/2}^N\}$ on a common probability space in the following way. Let $\mathfrak{T}$ be the set of sites in an infinite (rooted) 3-tree and consider a site percolation process on $\mathfrak{T}$ with each site being open, independently, with probability $1 - e^{-\beta p}$. We will call $C_0$ the corresponding percolation cluster containing 0. We also consider a collection $(B_i^N)_{i \in \mathfrak{T}, N>0}$ of independent Bernoulli random variables with $\mathbb{P}(B_i^N = 1) = \alpha_N$. With this, the random variable $\widetilde{\eta}_1^N(0)$, conditional on the event $\{0 \in H_N \cap \eta_{1/2}^N\}$, can be constructed as

$$\widetilde{\eta}_1^N(0) = \begin{cases} 1, & \text{if } B_i^N = 0 \text{ for all } i \in C_0 \cap B(0, (\log_2 N)/5), \\ 0, & \text{otherwise.} \end{cases}$$

It is clear that this construction gives the right conditional distribution for $\widetilde{\eta}_1^N(0)$.

Now let $l_N = \log_2(\alpha_N^{-1/2})$. Observe that $l_N < (\log_2 N)/5$ for large $N$, so we have that

$$\mathbb{P}(0 \in \widetilde{\eta}_1^N \mid 0 \in H_N \cap \eta_{1/2}^N)$$

$$\begin{aligned}
(3.2) \quad &= \mathbb{P}(0 \in \widetilde{\eta}_1^N, \operatorname{diam}(C_0) \leq l_N | 0 \in H_N \cap \eta_{1/2}^N) \\
&\quad + \mathbb{P}(0 \in \widetilde{\eta}_1^N, l_N < \operatorname{diam}(C_0) \leq \tfrac{1}{5} \log_2 N | 0 \in H_N \cap \eta_{1/2}^N) \\
&\quad + \mathbb{P}(0 \in \widetilde{\eta}_1^N, \operatorname{diam}(C_0) > \tfrac{1}{5} \log_2 N | 0 \in H_N \cap \eta_{1/2}^N).
\end{aligned}$$

For the first probability on the right-hand side we have that

$$\begin{aligned}
&\mathbb{P}(0 \in \widetilde{\eta}_1^N, \operatorname{diam}(C_0) \leq l_N | 0 \in H_N \cap \eta_{1/2}^N) \\
&\leq \mathbb{P}(0 < \operatorname{diam}(C_0) \leq l_N | 0 \in H_N \cap \eta_{1/2}^N) \\
&\longrightarrow \mathbb{P}(0 < \operatorname{diam}(C_0) < \infty | 0 \text{ is open}) = \frac{g_T(1 - e^{-\beta p})}{1 - e^{-\beta p}}.
\end{aligned}$$

This convergence is uniform in $p$ thanks to Lemma 3.3. On the other hand, since any subset of $\mathfrak{T}$ with diameter $n$ has at most $1 + 3 \cdot 2^{n-1} < 3 \cdot 2^n$ nodes, we get that

$$\begin{aligned}
&\mathbb{P}(0 \in \widetilde{\eta}_1^N, \operatorname{diam}(C_0) \leq l_N | 0 \in H_N \cap \eta_{1/2}^N) \\
&= \mathbb{P}(B_i^N = 0 \; \forall i \in C_0, \operatorname{diam}(C_0) \leq l_N | 0 \in H_N \cap \eta_{1/2}^N) \\
&= \mathbb{E}((1 - \alpha_N)^{|C_0|}, 0 < \operatorname{diam}(C_0) \leq l_N | 0 \in H_N \cap \eta_{1/2}^N) \\
&\geq (1 - \alpha_N)^{3\alpha_N^{-1/2}} \mathbb{P}(0 < \operatorname{diam}(C_0) \leq l_N | 0 \in H_N \cap \eta_{1/2}^N) \\
&\longrightarrow \frac{g_T(1 - e^{-\beta p})}{1 - e^{-\beta p}}
\end{aligned}$$

uniformly in $p$ by the same reason as above and because $(1 - \alpha_N)^{3\alpha_N^{-1/2}} \approx e^{-3\sqrt{\alpha_N}} \to 1$. We deduce that

$$\mathbb{P}(0 \in \widetilde{\eta}_1^N, \operatorname{diam}(C_0) \leq l_N | 0 \in H_N \cap \eta_{1/2}^N) \longrightarrow \frac{g_T(1 - e^{-\beta p})}{1 - e^{-\beta p}},$$

uniformly in $p$. For the second probability on the right-hand side of (3.2) we have that, since $\mathbb{P}(0 \in H_N \cap \eta_{1/2}^N) \geq C = (1 - e^{-\beta p})/2$ for large enough $N$,

$$\begin{aligned}
&\mathbb{P}(0 \in \widetilde{\eta}_1^N, l_N < \operatorname{diam}(C_0) \leq \tfrac{1}{5} \log_2 N | 0 \in H_N \cap \eta_{1/2}^N) \\
&\leq C^{-1} \mathbb{P}(l_N < \operatorname{diam}(C_0) \leq \tfrac{1}{5} \log_2 N) \\
&= C^{-1}[\mathbb{P}(\operatorname{diam}(C_0) > l_N) - \mathbb{P}(\operatorname{diam}(C_0) = \infty)] \\
&\quad - C^{-1}[\mathbb{P}(\operatorname{diam}(C_0) > \tfrac{1}{5} \log_2 N) - \mathbb{P}(\operatorname{diam}(C_0) = \infty)] \\
&\longrightarrow 0,
\end{aligned}$$

uniformly in $p$, again by Lemma 3.3. For the last probability in (3.2) we simply observe that

$$\mathbb{P}(0 \in \widetilde{\eta}_1^N, \operatorname{diam}(C_0) > \tfrac{1}{5} \log_2 N \,|\, 0 \in H_N \cap \eta_{1/2}^N)$$

$$\leq (1 - \alpha_N)^{1/5 \log_2 N} \approx e^{-1/5 \alpha_N \log_2 N} \longrightarrow 0.$$

The previous calculations and (3.2) imply that

$$\mathbb{P}(0 \in \widetilde{\eta}_1^N \,|\, 0 \in H_N \cap \eta_{1/2}^N) \longrightarrow \frac{g_T(1 - e^{-\beta p})}{1 - e^{-\beta p}},$$

uniformly in $p$. Putting this together with (3.1) we get the result. $\square$

Proof of Theorem 2. By Karr (1975), it is enough to prove that $\rho_0^N \Rightarrow p$ and that given any sequence $p_N$ in $[0, 1]$ converging to some $p' \in [0, 1]$, the sequence $\rho_1^N$, with $\eta_0^N$ started at a product measure of density $p_N$, converges weakly (or, equivalently, in probability) to $h_T(p')$.

The first part is straightforward. For the second part we will assume, for simplicity, that $p_N = p'$ for all $N$ and, moreover, that each site is occupied with probability $1 - e^{-\beta p'}$ after the growing season. The general case follows from the facts that $1 - (1 - \beta/N)^{p'N}$ converges uniformly as $N \to \infty$ to $1 - e^{-\beta p'}$ for $p' \in [0, 1]$ and that, by the preceding lemmas, all the convergences we will prove below are uniform on the initial density $p$.

Observe that by Markov's inequality, given any $\varepsilon > 0$

$$\mathbb{P}(|\rho_1^N - h_T(p')| > \varepsilon) \leq \frac{1}{\varepsilon} \mathbb{E}(|\rho_1^N - h_T(p')|),$$

so

$$
\begin{aligned}
\mathbb{P}(|\rho_1^N - h_T(p')| > \varepsilon) \leq {} & \frac{1}{\varepsilon} \mathbb{E}(|\rho_1^N - \widetilde{\rho}_1^N|) + \frac{1}{\varepsilon} \mathbb{E}\left( \left| \widetilde{\rho}_1^N - \frac{1}{N} |\widetilde{\eta}_1^N \cap H_N| \right| \right) \\
& + \frac{1}{\varepsilon} \mathbb{E}\left( \left| \frac{1}{N} |\widetilde{\eta}_1^N \cap H_N| - \mathbb{E}\left( \frac{1}{N} |\widetilde{\eta}_1^N \cap H_N| \right) \right| \right) \\
& + \frac{1}{\varepsilon} \left| \mathbb{E}\left( \frac{1}{N} |\widetilde{\eta}_1^N \cap H_N| \right) - h_T(p') \right|.
\end{aligned}
$$

(3.3)

Lemmas 3.1 and 3.4 imply that the first and last terms on the right-hand side of the inequality go to 0 as $N \to \infty$. The second one also goes to 0 since, using Lemma 3.2,

$$\mathbb{E}\left( \left| \widetilde{\rho}_1^N - \frac{1}{N} |\widetilde{\eta}_1^N \cap H_N| \right| \right) \leq \mathbb{E}\left( \frac{1}{N} |R_N \setminus H_N| \right) \longrightarrow 0.$$

To deal with the third term, observe that

$$\operatorname{Var}(|\widetilde{\eta}_1^N \cap H_N|)$$

$$= \sum_{i=0}^{N-1} \sum_{j=0}^{N-1} \mathrm{Cov}(\mathbf{1}_{i \in \widetilde{\eta}_1^N \cap H_N}, \mathbf{1}_{j \in \widetilde{\eta}_1^N \cap H_N})$$

$$\leq |\{(i,j) \in H_N \times H_N : B(i, (\log_2 N)/5) \cap B(j, (\log_2 N)/5) \neq \varnothing\}|$$

$$= |\{(i,j) \in H_N \times H_N : |i-j| \leq 2(\log_2 N)/5\}| \leq 2N \cdot N^{2/5}.$$

Hence, by Jensen's inequality,

$$\mathbb{E}\left(\left|\frac{1}{N}|\widetilde{\eta}_1^N \cap H_N| - \mathbb{E}\left(\frac{1}{N}|\widetilde{\eta}_1^N \cap H_N|\right)\right|\right)^2 \leq \mathrm{Var}\left(\frac{1}{N}|\widetilde{\eta}_1^N \cap H_N|\right)$$

$$\leq \frac{2N \cdot N^{2/5}}{N^2} \longrightarrow 0.$$

We deduce from (3.3) that $\rho_1^N$ converges in probability to $h_T(p')$. □

**4. Proof of Theorem 3.** As in the case of the 3-tree, we let $a_0$ be the solution of $f(a_0) = p_c$ [i.e., $a_0 = \log(1/(1-p_c))/\beta$] and $a_1 = h_L(p_c)$ (see Figure 1 for a sketch of these values in the case of the 3-tree). It is enough to prove, by Lasota and Yorke (1973), that there is a $K \in \mathbb{N}$ such that

$$(4.1) \qquad \inf_{p \in [a_1, p_c] \setminus \{a_0\}} |(h_T^K)'(p)| > 1.$$

Fix any $\beta_1 > \beta_c$. Since $a_1$ is bounded away from 0 for $\beta \in (\beta_c, \beta_1)$, there is a $K \in \mathbb{N}$ such that $\min\{k \in \mathbb{N} : f^k(a_1) > p_c\} \leq K - 1$ for any such $\beta$. In particular, since $a_0$ is always less than $p_c$ we deduce that given any $\beta \in (\beta_c, \beta_1)$ and any $p \in [a_1, p_c]$, the $K$-tuple $(p, h_L(p), \ldots, h_L^{K-1}(p))$ contains at least one point in $(a_0, p_c]$.

Now recall that $f'' < 0$, so $f'$ attains its minimum on the interval $[a_1, a_0]$ at $a_0$, and at this point its value is $\beta(1-p_c)$. Thus for every $\beta \in (\beta_c, \beta_1)$, this minimum is larger than $\beta_c(1-p_c)$. Since $g_L(p) = p$ for $p \in [a_1, a_0]$ we deduce that

$$|h_L(p)| \geq \beta_c(1-p_c) \qquad \text{for all } p \in [a_1, a_0].$$

Now using the fact that $a_0 \uparrow p_c$ as $\beta \downarrow \beta_c$, we can choose given any $\varepsilon > 0$ a $\beta_2 \in (\beta_c, \beta_1)$ so that $f(p_c) - p_c = f(p_c) - f(a_0) < \varepsilon$ for any $\beta \in (\beta_c, \beta_2)$. Since (1.2) implies that

$$g_L'(p) = 1 - \theta_L'(p) \xrightarrow[p \downarrow p_c]{} -\infty,$$

we can choose a small enough $\varepsilon$, so that

$$|h_L'(p)| = |g_L'(f(p))||f(p)| > \max\{[\beta_c(1-p_c)]^{-(K-1)}, 1\}$$

for all $p \in (a_0, p_c]$, and thus this inequality holds for all $\beta \in (\beta_c, \beta_2)$.

Putting the previous arguments together with the fact that

$$(h_L^K)'(p) = h_L'(h_L^{K-1}(p))h_L'(h_L^{K-2}(p)) \cdots h_L'(p)$$

we deduce that (4.1) holds for all $\beta \in (\beta_c, \beta_2)$.

**5. Proof of Theorem 4.**  Given $i \in \mathbb{T}_N$ and $m \in \mathbb{N}$ we will write

$$B(i,m) = \{j \in \mathbb{T}_N : \|i - j\|_\infty \leq m\} \quad \text{and} \quad V(m) = (2m+1)^d = |B(i,m)|$$

(here and in what follows all differences $i - j$ for $i, j \in \mathbb{T}_N$ are computed modulo $N$). Define, for $k \in \mathbb{N}$,

$$d_k^N(i) = \frac{1}{V(r_N)} \sum_{\|j - i\|_\infty \leq r_N} \eta_k^N(j)$$

and

$$G_k^N(\varepsilon) = \{i \in \mathbb{T}_N : |d_k^N(i) - h_L^k(p)| < \varepsilon\}.$$

$d_k^N(i)$ is the density of occupied sites in the growth neighborhood of $i$, while $G_k^N(\varepsilon)$ can be thought of as the set of "good sites at time $k$," where a site is said to be good at time $k$ if the density of occupied sites in its growth neighborhood at that time is close to the desired value $h_L^k(p)$. The proof of Theorem 4 will depend on the following proposition:

PROPOSITION 5.1.  *Fix $\varepsilon_1, \varepsilon_2 > 0$ and $k \in \mathbb{N}$ and assume that*

$$(5.1) \qquad\qquad \frac{1}{N^d} \mathbb{E}(|\mathbb{T}_N \setminus G_k^N(\delta_1)|) < \delta_2.$$

*Then if $\delta_1$ and $\delta_2$ are small enough and $N$ is large enough,*

$$\frac{1}{N^d} \mathbb{E}(|\mathbb{T}_N \setminus G_{k+1}^N(\varepsilon_1)|) < \varepsilon_2.$$

This result will allow us to give an inductive proof of Theorem 4. We will need thus the following lemma:

LEMMA 5.2.  *Given any $\delta > 0$,*

$$\frac{1}{N^d} \mathbb{E}(|\mathbb{T}_N \setminus G_0^N(\delta)|) \xrightarrow[N \to \infty]{} 0.$$

PROOF.  By translation invariance,

$$\mathbb{E}(|\mathbb{T}_N \setminus G_0^N(\delta)|) = \sum_{i \in \mathbb{T}_N} \mathbb{P}(i \notin G_0^N(\delta)) = N^d \mathbb{P}(|d_0^N(0) - p| \geq \delta).$$

Since $\mathbb{E}(d_0^N(0)) = p$, Chebyshev's inequality and the fact that (by definition) $V(r_N)d_0^N(0)$ is the sum of $V(r_N)$ independent Bernoulli random variables with success probability $p$ imply that

$$\mathbb{P}(|d_0^N(0) - p| \geq \delta) \leq \frac{1}{\delta^2 V(r_N)^2} V(r_N)p(1-p),$$

so

$$\frac{1}{N^d}\mathbb{E}(|\mathbb{T}_N \setminus G_0^N(\delta)|) \leq \frac{1}{\delta^2 V(r_N)} p(1-p) \longrightarrow 0. \qquad \square$$

Now we turn to the proof of Proposition 5.1. Many parts in the argument will be similar to those in the proof of Theorem 2 and the lemmas that preceded it, so we will skip some details. We begin with some preliminary results. Throughout this part, and until the proof of Theorem 4, we fix $k, \delta_1, \delta_2, \varepsilon_1, \varepsilon_2$ and assume that (5.1) holds.

Observe that since each occupied site $i$ sends a Poisson[$\beta$] number of births during the growing season, each to a site chosen randomly from $B(i, r_N)$, we can equivalently think of each occupied site $i$ as sending a Poisson[$\beta/V(r_N)$] number of births to each of its $V(r_N)$ neighbors at a distance smaller than $r_N$. Hence during the growing season, each site $i$ receives $\sum_{\|j-i\|_\infty \leq r_N} \eta_k^N(j) Y_{j,i}$ births, where $(Y_{i,j})_{i,j \in \mathbb{T}_N}$ are i.i.d. Poisson[$\beta/V(r_N)$] random variables. Conditional on $d_k^N(i)$, this last sum is distributed as a Poisson[$d_k^N(i)\beta$] random variable. We deduce that we can regard the growing season as taking place as follows:

Given $\eta_k^N$, each $i$ will be in $\eta_{k+1/2}^N$ with probability equal to the probability that a Poisson[$d_k^N(i)\beta$] random variable is positive, that is, with probability $1 - e^{-\beta d_k^N(i)}$.

The Poisson random variables above are taken to be independent of each other.

Let $l_N = \sqrt{r_N/\alpha_N}$ and observe that

$$\frac{l_N}{r_N} = \frac{1}{\sqrt{\alpha_N r_N}} \longrightarrow 0 \quad \text{and} \quad \alpha_N l_N = \sqrt{\alpha_N r_N} \longrightarrow \infty.$$

We let $\widetilde{\eta}_{k+1}^N$ be the configuration obtained from $\eta_{k+1/2}^N$ by ignoring infections coming from a distance greater than $l_N$.

LEMMA 5.3.

$$\frac{1}{N^d} \sum_{i \in \mathbb{T}_N} \mathbb{E}(|\eta_{k+1}^N(i) - \widetilde{\eta}_{k+1}^N(i)|) \xrightarrow[N \to \infty]{} 0.$$

In particular,

$$\mathbb{E}(|\rho_{k+1}^N - \widetilde{\rho}_{k+1}^N|) \longrightarrow 0.$$

PROOF.   By translation invariance, and repeating the arguments of the proof of Lemma 3.1, we get that

$$\frac{1}{N^d} \sum_{i \in \mathbb{T}_N} \mathbb{E}(|\eta_{k+1}^N(i) - \widetilde{\eta}_{k+1}^N(i)|) = \frac{1}{N^d} \sum_{i \in \mathbb{T}_N} \mathbb{P}(\eta_{k+1}^N(i) \neq \widetilde{\eta}_{k+1}^N(i))$$

$$= \mathbb{P}(0 \in \widetilde{\eta}_{k+1}^N \setminus \eta_{k+1}^N) \leq (1 - \alpha_N)^{l_N}$$

$$\approx e^{-\alpha_N l_N} \longrightarrow 0. \qquad \square$$

Before continuing, it is useful to give an explicit construction of the random variable $\widetilde{\eta}_{k+1}^N(0)$. Consider a collection $X = (X(i))_{i \in \mathbb{Z}^d}$ of i.i.d. random variables with uniform distribution in $[0, 1]$ and, given $\eta_k^N$, construct $\eta_{k+1/2}^N$ as follows:

$$\eta_{k+1/2}^N(i) = \mathbf{1}_{X(i) > e^{-\beta d_k^N(i)}}.$$

Observe that with this choice, $\mathbb{P}(\eta_{k+1/2}^N(i) = 1) = 1 - e^{-\beta d_k^N(i)}$ as required. We will call $C_0^N$ the open cluster in $\eta_{k+1/2}^N$ containing 0. Define $(B_i^N)_{i \in \mathbb{Z}^d, N>0}$ as in Section 3 and set

$$\widetilde{\eta}_{k+1}^N(0) = \begin{cases} 1, & \text{if } \widetilde{\eta}_{k+1/2}^N(0) = 1 \text{ and } B_i^N = 0 \text{ for all } i \in C_0^N \cap B(0, l_N), \\ 0, & \text{otherwise.} \end{cases}$$

This construction gives the right distribution for $\widetilde{\eta}_{k+1}^N(0)$.

We introduce another modification of $\eta_{k+1}^N$: let $\widehat{\eta}_{k+1}^N$ be the configuration obtained from $\eta_k^N$ in the same way as $\widetilde{\eta}_{k+1}^N$, except that in the growing season we replace $\eta_{k+1/2}^N$ by the configuration $\widehat{\eta}_{k+1/2}^N$ defined by

$$\widehat{\eta}_{k+1/2}^N(i) = \mathbf{1}_{X(i) > e^{-\beta h_L^k(p)}}$$

(using the same family of variables $X$). That is, $\widehat{\eta}_{k+1/2}^N$ corresponds to running the growth step as if the density of occupied sites in the ball of radius $r_N$ around each site was exactly $h_L^k(p)$. $\widehat{\rho}_k^N$ will denote the density of occupied sites in this modified process, that is, $\widehat{\rho}_k^N = |\widehat{\eta}_k^N|/N^d$. We will call $C_0$ the open cluster containing 0 in the site percolation process in all of $\mathbb{Z}^d$ constructed from the collection of random variables $X$ with each site being open with probability $1 - e^{-\beta h_L^k(p)}$.

LEMMA 5.4.   *Given any $\varepsilon > 0$, if $\delta_1$ and $\delta_2$ are small enough, then*

$$\mathbb{E}(|\widetilde{\rho}_{k+1}^N - \widehat{\rho}_{k+1}^N|) \leq \varepsilon.$$

PROOF.   The idea behind the proof of this result is the following. By (5.1), the density of occupied sites is close to $h_L^k(p)$ around most sites. If this holds for some site $i$, then in a box around $i$ the density must still be close to this. We then prove the result by comparing $\widetilde{\eta}_{k+1}^N$ and $\widehat{\eta}_{k+1}^N$ with processes in which the outcome of the growth step is replaced by product measures of sligthly smaller and slightly larger densities.

To get started we observe that

$$\mathbb{E}(|\widetilde{\rho}_{k+1}^N - \widehat{\rho}_{k+1}^N|)$$

$$\leq \frac{1}{N^d} \sum_{i \in \mathbb{T}_N} \mathbb{E}(|\widetilde{\eta}_{k+1}^N(i) - \widehat{\eta}_{k+1}^N(i)|) = \mathbb{P}(\widetilde{\eta}_{k+1}^N(0) \neq \widehat{\eta}_{k+1}^N(0))$$

(5.2)

$$\leq \mathbb{P}(\widetilde{\eta}_{k+1}^N(0) \neq \widehat{\eta}_{k+1}^N(0), 0 \in G_k^N(\delta_1)) + \mathbb{P}(0 \notin G_k^N(\delta_1))$$

$$\leq \mathbb{P}(\widetilde{\eta}_{k+1}^N(0) \neq \widehat{\eta}_{k+1}^N(0), 0 \in G_k^N(\delta_1)) + \delta_2,$$

where in last bound we used (5.1). To deal with the last probability we first observe that given any $i \in B(0, l_N)$,

$$d_k^N(i) = \frac{1}{V(r_N)} \sum_{j \in B(i, r_N)} \eta_k^N(j)$$

$$= d_k^N(0) + \frac{1}{V(r_N)} \sum_{j \in B(i, r_N) \setminus B(0, r_N)} \eta_k^N(j)$$

$$- \frac{1}{V(r_N)} \sum_{j \in B(0, r_N) \setminus B(i, r_N)} \eta_k^N(j)$$

$$\leq d_k^N(0) + \frac{|B(i, r_N) \setminus B(0, r_N)|}{V(r_N)}$$

and thus, since the cardinality in the last term is largest when $i$ is at any of the $2^d$ corners of the hypercube $B(0, l_N)$, we have that for some $C > 0$

$$|d_k^N(i) - d_k^N(0)| \leq C \frac{r_N^{d-1} l_N}{V(r_N)} \approx \frac{l_N}{r_N} \longrightarrow 0.$$

We deduce that

$$\mathbb{P}(\widetilde{\eta}_{k+1}^N(0) \neq \widehat{\eta}_{k+1}^N(0), 0 \in G_k^N(\delta_1))$$

$$\leq \mathbb{P}(\widetilde{\eta}_{k+1}^N(0) \neq \widehat{\eta}_{k+1}^N(0), |d_k^N(i) - h_L^k(p)| \leq 2\delta_1 \ \forall i \in B(0, l_N), 0 \in G_k^N(\delta_1))$$

$$+ \mathbb{P}(|d_k^N(i) - h_L^k(p)| > 2\delta_1 \text{ for some } i \in B(0, l_N), 0 \in G_k^N(\delta_1))$$

$$\leq \mathbb{P}(\widetilde{\eta}_{k+1}^N(0) \neq \widehat{\eta}_{k+1}^N(0), |d_k^N(i) - h_L^k(p)| \leq 2\delta_1 \ \forall i \in B(0, l_N))$$

$$+ \mathbb{P}(|d_k^N(i) - d_k^N(0)| > \delta_1 \text{ for some } i \in B(0, l_N))$$

$$+ \mathbb{P}(|d_k^N(0) - h_L^k(p)| > \delta_1, 0 \in G_k^N(\delta_1))$$

$$= \mathbb{P}(\widetilde{\eta}_{k+1}^N(0) \neq \widehat{\eta}_{k+1}^N(0), |d_k^N(i) - h_L^k(p)| \leq 2\delta_1 \ \forall i \in B(0, l_N))$$

for large enough $N$.

Next, we introduce the following notation: $\xi_{1/2}^q$ will be the set of open sites in a site percolation process in $\mathbb{Z}^d$ with each site being open with probability $1 - e^{-\beta q}$ for $q \in [0,1]$ constructed from the family of random variables $X$. In other words, we put $\xi_{1/2}^q(i) = \mathbf{1}_{X(i) > e^{-\beta q}}$ for each $i \in \mathbb{Z}^d$. We also let $\xi_1^{q,N} \subseteq \mathbb{T}_N$ be the configuration obtained after running the epidemic step on $\xi_{1/2}^q \cap \mathbb{T}_N$ (this is done on the torus $\mathbb{T}_N$, so we take into account the periodic boundary conditions of the torus while running the epidemic), using the variables $(B_i^N)_{i \in \mathbb{T}_N}$, and ignoring infections coming from a distance greater than $l_N$. Observe that with these definitions, $\widehat{\eta}_{k+1/2}^N = \xi_{1/2}^{h_L^k(p)} \cap \mathbb{T}_N$ and $\widehat{\eta}_{k+1}^N = \xi_1^{h_L^k(p),N}$. The key fact is the following:

$$\mathbb{P}(\widetilde{\eta}_{k+1}^N(0) \neq \widehat{\eta}_{k+1}^N(0), |d_k^N(i) - h_L^k(p)| \leq 2\delta_1 \ \forall i \in B(0, l_N))$$

$$\leq \mathbb{P}(\xi_1^{h_L^k(p)+2\delta_1, N}(0) = 0, \xi_1^{h_L^k(p),N}(0) = 1)$$

$$\text{(5.3)} \qquad + \mathbb{P}(\xi_{1/2}^{h_L^k(p)-2\delta_1}(0) = 0, \xi_{1/2}^{h_L^k(p)}(0) = 1)$$

$$+ \mathbb{P}(\xi_1^{h_L^k(p)-2\delta_1, N}(0) = 1, \xi_1^{h_L^k(p),N}(0) = 0)$$

$$+ \mathbb{P}(\xi_{1/2}^{h_L^k(p)+2\delta_1}(0) = 1, \xi_{1/2}^{h_L^k(p)-2\delta_1}(0) = 0).$$

To see that this is true observe that $|d_k^N(i) - h_L^k(p)| \leq 2\delta_1$ for all $i \in B(0, l_N)$ implies that

$$1 - e^{-\beta(h_L^k(p)-2\delta_1)} \leq 1 - e^{-\beta d_k^N(i)} \leq 1 - e^{-\beta(h_L^k(p)+2\delta_1)}$$

for all $i \in B(0, l_N)$, and thus

$$\xi_{1/2}^{h_L^k(p)-2\delta_1} \cap B(0, l_N) \subseteq C_0^N \cap B(0, l_N) \subseteq \xi_{1/2}^{h_L^k(p)+2\delta_1} \cap B(0, l_N).$$

Assuming this, we have that $\widetilde{\eta}_{k+1}^N(0) = 0$ and $\widehat{\eta}_{k+1}^N(0) = 1$ implies that $\xi_1^{h_L^k(p),N}(0) = \xi_{1/2}^{h_L^k(p)}(0) = 1$, and either $\widetilde{\eta}_{k+1/2}^N(0) = 0$, which implies that $\xi_{1/2}^{h_L^k(p)-2\delta_1}(0) = 0$, or $\widetilde{\eta}_{k+1/2}^N(0) = 1$ but there is an infection in $C_0^N \cap B(0, l_N)$, which implies that $\xi_1^{h_L^k(p)+2\delta_1, N}(0) = 0$. Similarly, $\widetilde{\eta}_{k+1}^N(0) = 1$ and $\widehat{\eta}_{k+1}^N = 0$ implies that $\xi_{1/2}^{h_L^k(p)+2\delta_1}(0) = 1$, $\xi_1^{h_L^k(p),N} = 0$, and there is no infection in $C_0^N \cap B(0, l_N)$, and thus $\xi_1^{h_L^k(p)-2\delta_1, N}(0) = 1$ whenever $\xi_{1/2}^{h_L^k(p)-2\delta_1}(0) = 1$.

To finish the proof we need to bound the probabilities on the right-hand side of (5.3). For the first one, since $\xi_{1/2}^{h_L^k(p)-2\delta_1} \subseteq \xi_{1/2}^{h_L^k(p)} \subseteq \xi_{1/2}^{h_L^k(p)+2\delta_1}$, we have that if $\#\xi$ denotes the size of the cluster containing 0 in the configuration given by $\xi$, then

$$\mathbb{P}(\xi_1^{h_L^k(p)+2\delta_1, N}(0) = 0, \xi_1^{h_L^k(p),N}(0) = 1)$$

$$\leq \mathbb{P}(\xi_1^{h_L^k(p)+2\delta_1,N}(0)=0, \xi_1^{h_L^k(p),N}(0)=1, \#\xi_{1/2}^{h_L^k(p)+2\delta_1} < \infty)$$

$$+ \mathbb{P}(\xi_1^{h_L^k(p),N}(0)=1, \#\xi_{1/2}^{h_L^k(p)}=\infty) + \mathbb{P}(\#\xi_{1/2}^{h_L^k(p)} < \#\xi_{1/2}^{h_L^k(p)+2\delta_1}=\infty).$$

The first probability on the right-hand side is bounded by

$$(5.4) \qquad \mathbb{P}(\xi_1^{h_L^k(p)+2\delta_1,N}(0)=0, \xi_{1/2}^{h_L^k(p)+2\delta_1}(0)=1, \#\xi_{1/2}^{h_L^k(p)+2\delta_1} < \infty)$$

$$\leq \mathbb{E}(1-(1-\alpha_N)^{\#\xi_{1/2}^{h_L^k(p)+2\delta_1}}, \#\xi_{1/2}^{h_L^k(p)+2\delta_1} < \infty),$$

which goes to 0 by the dominated convergence theorem. The second one goes to 0 as well because it is bounded by $(1-\alpha_N)^{l_N} \approx e^{-\alpha_N l_N}$. The third one equals

$$\theta_L(h_L^k(p) + 2\delta_1) - \theta_L(h_L^k(p)),$$

which is less than $\varepsilon/2$ for small enough $\delta_1$ by the (uniform) continuity of the percolation probability $\theta_L(p)$ for $p \in [0,1]$. The other two probabilities on the right-hand side of (5.3) can be bounded similarly, yielding

$$\mathbb{P}(\widetilde{\eta}_{k+1}^N(0) \neq \widehat{\eta}_{k+1}^N(0), 0 \in G_k^N(\delta_1)) < \varepsilon$$

for large enough $N$ and small enough $\delta_1$. Putting this together with (5.2) gives the result. $\square$

LEMMA 5.5.

$$|\mathbb{E}(\widehat{\rho}_{k+1}^N) - h_L^{k+1}(p)| \longrightarrow 0.$$

PROOF. This proof is similar to that of Lemma 3.4. First we observe that

$$(5.5) \qquad \mathbb{E}(\widehat{\rho}_{k+1}^N) = \mathbb{P}(0 \in \widehat{\eta}_{k+1}^N | 0 \in \widehat{\eta}_{k+1/2}^N) \mathbb{P}(0 \in \widehat{\eta}_{k+1/2}^N)$$

$$= \mathbb{P}(0 \in \widehat{\eta}_{k+1}^N | 0 \in \widehat{\eta}_{k+1/2}^N)[1 - e^{-\beta h_L^k(p)}]$$

and

$$(5.6) \quad \mathbb{P}(0 \in \widehat{\eta}_{k+1}^N, \mathrm{diam}(C_0)=\infty | 0 \in \widehat{\eta}_{k+1/2}^N) \leq (1-\alpha_N)^{l_N} \approx e^{-\alpha_N l_N} \longrightarrow 0.$$

Now

$$(5.7) \qquad \mathbb{P}(0 \in \widehat{\eta}_{k+1}^N, \mathrm{diam}(C_0) < \infty | 0 \in \widehat{\eta}_{k+1/2}^N)$$

$$= \mathbb{P}(0 \in \widehat{\eta}_{k+1}^N, \mathrm{diam}(C_0) \leq l_N | 0 \in \widehat{\eta}_{k+1/2}^N)$$

$$+ \mathbb{P}(0 \in \widehat{\eta}_{k+1}^N, l_N < \mathrm{diam}(C_0) < \infty | 0 \in \widehat{\eta}_{k+1/2}^N)$$

and, trivially,

$$(5.8) \qquad \mathbb{P}(0 \in \widehat{\eta}^N_{k+1}, l_N < \mathrm{diam}(C_0) < \infty | 0 \in \widehat{\eta}^N_{k+1/2})$$
$$\leq \mathbb{P}(l_N < \mathrm{diam}(C_0) < \infty | 0 \in \widehat{\eta}^N_{k+1/2}) \longrightarrow 0.$$

On the other hand,

$$\mathbb{P}(0 \in \widehat{\eta}^N_{k+1}, \mathrm{diam}(C_0) \leq l_N | 0 \in \widehat{\eta}^N_{k+1/2})$$
$$= \mathbb{P}(B^N_i = 0 \ \forall i \in C_0 \cap B(0, l_N), \mathrm{diam}(C_0) \leq l_N | 0 \text{ is open})$$
$$= \mathbb{E}((1 - \alpha_N)^{|C_0 \cap B(0, l_N)|}, \mathrm{diam}(C_0) \leq l_N | 0 \text{ is open})$$
$$= \mathbb{P}(\mathrm{diam}(C_0) \leq l_N | 0 \text{ is open})$$
$$\quad - \mathbb{E}(1 - (1 - \alpha_N)^{|C_0 \cap B(0, l_N)|}, \mathrm{diam}(C_0) \leq l_N | 0 \text{ is open}).$$

The second expectation is positive and bounded from above by

$$\mathbb{E}(1 - (1 - \alpha_N)^{|C_0|}, |C_0| < \infty | 0 \text{ is open}),$$

so it goes to 0 as $N \to \infty$ by the dominated convergence theorem as in (5.4). Thus

$$\lim_{N \to \infty} \mathbb{P}(0 \in \widehat{\eta}^N_{k+1}, \mathrm{diam}(C_0) \leq l_N | 0 \in \widehat{\eta}^N_{k+1/2})$$
$$= \mathbb{P}(\mathrm{diam}(C_0) < \infty | 0 \text{ is open})$$
$$= \frac{\mathbb{P}(0 < \mathrm{diam}(C_0) < \infty)}{1 - e^{-\beta h^k_L(p)}} = \frac{g_L(1 - e^{-\beta h^k_L(p)})}{1 - e^{-\beta h^k_L(p)}}.$$

Putting this together with (5.7) and (5.8) we get that

$$\left| \mathbb{P}(0 \in \widehat{\eta}^N_{k+1}, \mathrm{diam}(C_0) < \infty | 0 \in \widehat{\eta}^N_{k+1/2}) - \frac{h^{k+1}_L(p)}{1 - e^{-\beta h^k_L(p)}} \right| \longrightarrow 0$$

and thus by (5.5) and (5.6) we obtain

$$|\mathbb{E}(\widehat{\rho}^N_{k+1}) - h^{k+1}_L(p)| \longrightarrow 0$$

as required. □

PROOF OF PROPOSITION 5.1.

$$\frac{1}{N^d} \mathbb{E}(|\mathbb{T}_N \setminus G^N_{k+1}(\varepsilon_1)|) = \mathbb{P}(0 \notin G^N_{k+1}(\varepsilon_1)) = \mathbb{P}(|d^N_{k+1}(0) - h^{k+1}_L(p)| \geq \varepsilon_1)$$
$$\leq \frac{1}{\varepsilon_1} \mathbb{E}(|d^N_{k+1}(0) - h^{k+1}_L(p)|).$$

Hence

$$\frac{1}{N^d}\mathbb{E}(|\mathbb{T}_N \setminus G_{k+1}^N(\varepsilon_1)|)$$

(5.9)
$$\leq \frac{1}{\varepsilon_1}[\mathbb{E}(|d_{k+1}^N(0) - \widetilde{d}_{k+1}^N(0)|) + \mathbb{E}(|\widetilde{d}_{k+1}^N(0) - \mathbb{E}(\widetilde{\rho}_{k+1}^N)|)$$
$$+ \mathbb{E}(|\mathbb{E}(\widetilde{\rho}_{k+1}^N) - \widetilde{\rho}_{k+1}^N|) + \mathbb{E}(|\widetilde{\rho}_{k+1}^N - h_L^{k+1}(p)|)],$$

where $\widetilde{d}_{k+1}^N(0) = \frac{1}{V(r_N)}\sum_{\|j\|_\infty \leq r_N}\widetilde{\eta}_{k+1}^N(j)$.

For fixed $\varepsilon > 0$ we want to show that each of the expectations on the right-hand side of the last inequality can be bounded by $\varepsilon$ if $N$ is large enough and $\delta_1$ and $\delta_2$ are small enough. The bound for the last one follows directly from the triangle inequality and Lemmas 5.4 and 5.5.

For the first one we have by translation invariance that

$$\mathbb{E}(|d_{k+1}^N(0) - \widetilde{d}_{k+1}^N(0)|) \leq \frac{1}{V(r_N)}\sum_{\|j\|_\infty \leq r_N}\mathbb{E}(|\eta_{k+1}^N(j) - \widetilde{\eta}_{k+1}^N(j)|)$$

$$= \frac{1}{N^d V(r_N)}\sum_{i \in \mathbb{T}_N}\sum_{j \in B(i,r_N)}\mathbb{E}(|\eta_{k+1}^N(j) - \widetilde{\eta}_{k+1}^N(j)|)$$

$$= \frac{1}{N^d V(r_N)}\sum_{j \in \mathbb{T}_N}\mathbb{E}\left(\sum_{i \in B(j,r_N)}|\eta_{k+1}^N(j) - \widetilde{\eta}_{k+1}^N(j)|\right)$$

$$= \frac{1}{N^d}\sum_{j \in \mathbb{T}_N}\mathbb{E}(|\eta_{k+1}^N(j) - \widetilde{\eta}_{k+1}^N(j)|) < \varepsilon$$

for large enough $N$ by Lemma 5.3.

For the second one we first observe that, again by translation invariance, $\mathbb{E}(\widetilde{d}_{k+1}^N(0)) = \mathbb{E}(\widetilde{\rho}_{k+1}^N)$. Hence

(5.10)
$$\mathbb{E}(|\widetilde{d}_{k+1}^N(0) - \mathbb{E}(\widetilde{\rho}_{k+1}^N)|)^2$$
$$\leq \mathrm{Var}(\widetilde{d}_{k+1}^N(0))$$
$$= \frac{1}{V(r_N)^2}\sum_{i,j \in B(0,r_N)}\mathrm{Cov}(\widetilde{\eta}_{k+1}^N(i), \widetilde{\eta}_{k+1}^N(j))$$
$$\leq \frac{1}{V(r_N)^2}|\{i, j \in B(0,r_N): \|i-j\|_\infty \leq l_N\}| \approx \frac{V(l_N)}{V(r_N)} \longrightarrow 0.$$

The bound for the third expectation on the right-hand side of (5.9) follows from the exact same argument as previous one. We deduce that

$$\frac{1}{N^d}\mathbb{E}(|\mathbb{T}_N \setminus G_{k+1}^N(\varepsilon_1)|) \leq \frac{4\varepsilon}{\varepsilon_1}$$

for large enough $N$, and thus choosing $\varepsilon < \varepsilon_1 \varepsilon_2 / 4$ gives the result.    □

PROOF OF THEOREM 4.  Since $[0, 1]$ is compact, it is enough to prove the convergence of the finite-dimensional distributions of $\rho_k^N$, and since our limit is deterministic, we only need to prove that

$$(5.11) \qquad \mathbb{P}(|\rho_k^N - h_L^k(p)| > \varepsilon) \xrightarrow[N \to \infty]{} 0$$

for every $k \geq 0$ and $\varepsilon > 0$. Proceeding as in the proof of Theorem 2 we have that

$$(5.12) \begin{aligned} \mathbb{P}(|\rho_k^N - h_L^k(p)| > \varepsilon) &\leq \frac{1}{\varepsilon} \mathbb{E}(|\rho_k^N - \widetilde{\rho}_k^N|) + \frac{1}{\varepsilon} \mathbb{E}(|\widetilde{\rho}_k^N - \widehat{\rho}_k^N|) \\ &\quad + \frac{1}{\varepsilon} \mathbb{E}(|\widehat{\rho}_k^N - \mathbb{E}(\widehat{\rho}_k^N)|) + \frac{1}{\varepsilon} |\mathbb{E}(\widehat{\rho}_k^N) - h_L(p)|. \end{aligned}$$

By Lemmas 5.3, 5.4 and 5.5, given any $\upsilon > 0$ there are constants $\delta_1^{k-1}, \delta_2^{k-1} > 0$ such that

$$(5.13) \qquad \frac{V(l_N)}{N^d} \mathbb{E}(|\mathbb{T}_N \setminus G_{k-1}^N(\delta_1^{k-1})|) < \delta_2^{k-1}$$

implies that the first, second and last terms on the right-hand side of (5.12) are each bounded by $\upsilon \varepsilon$ for large enough $N$. The third term is also less than $\upsilon \varepsilon$ for large $N$, which follows from repeating again the argument in (5.10). We deduce that

$$(5.14) \qquad \mathbb{P}(|\rho_k^N - h_L^k(p)| > \varepsilon) < 4\upsilon$$

for large enough $N$ provided that (5.13) holds.

Similarly, Proposition 5.1 implies that (5.13) will hold provided that

$$\frac{V(l_N)}{N^d} \mathbb{E}(|\mathbb{T}_N \setminus G_{k-2}^N(\delta_1^{k-2})|) < \delta_2^{k-2}$$

for some $\delta_1^{k-2}, \delta_2^{k-2} > 0$. Repeating this procedure inductively we deduce that (5.14) holds provided that

$$\frac{V(l_N)}{N^d} \mathbb{E}(|\mathbb{T}_N \setminus G_0^N(\delta_1^0)|) < \delta_2^0$$

for some small $\delta_1^0, \delta_2^0 > 0$, which holds for large enough $N$ by Lemma 5.2, and thus (5.11) follows.    □

## REFERENCES

DURRETT, R. (1995). Ten lectures on particle systems. In *Lectures on Probability Theory (Saint–Flour, 1993). Lecture Notes in Mathematics* **1608** 97–201. Springer, Berlin. MR1383122

DURRETT, R. (2007). *Random Graph Dynamics*. Cambridge Univ. Press, Cambridge. MR2271734

GRIMMETT, G. (1999). *Percolation*, 2nd ed. *Grundlehren der Mathematischen Wissenschaften [Fundamental Principles of Mathematical Sciences]* **321**. Springer, Berlin. MR1707339

JANSON, S., LUCZAK, T. and RUCINSKI, A. (2000). *Random Graphs*. Wiley-Interscience, New York. MR1782847

KARR, A. F. (1975). Weak convergence of a sequence of Markov chains. *Z. Wahrsch. Verw. Gebiete* **33** 41–48. MR0394882

KESTEN, H. and ZHANG, Y. (1987). Strict inequalities for some critical exponents in two-dimensional percolation. *J. Statist. Phys.* **46** 1031–1055. MR893131

LASOTA, A. and YORKE, J. A. (1973). On the existence of invariant measures for piecewise monotonic transformations. *Trans. Amer. Math. Soc.* **186** 481–488. MR0335758

LI, T. Y. and YORKE, J. A. (1975). Period three implies chaos. *Amer. Math. Monthly* **82** 985–992. MR0385028

RUDIN, W. (1976). *Principles of Mathematical Analysis*, 3rd ed. McGraw-Hill, New York. MR0385023

DEPARTMENT OF MATHEMATICS
CORNELL UNIVERSITY
523 MALOTT HALL
ITHACA, NEW YORK 14853
USA
E-MAIL: rtd1@cornell.edu

CENTER FOR APPLIED MATHEMATICS
CORNELL UNIVERSITY
657 RHODES HALL
ITHACA, NEW YORK 14853
USA
E-MAIL: dir4@cornell.edu